\documentclass[11pt,a4paper]{article}
\usepackage[T1]{fontenc}
\usepackage{lmodern}
\usepackage[margin=27mm]{geometry}
\usepackage{amsmath,amssymb,amsthm,mathtools}
\usepackage{microtype}
\usepackage{enumitem}
\usepackage{xcolor}
\usepackage{hyperref}
\hypersetup{colorlinks=true,linkcolor=blue!45!black,citecolor=blue!45!black,urlcolor=blue!45!black,
 pdfauthor={Long Tian, Yongshun Yang},
 pdftitle={Sharp box-dimension bounds for fractional integrals under unequal Holder and anti-Holder exponents},
 pdfsubject={Riemann--Liouville, Hadamard and Weyl integrals; sharp graph box-dimension bounds}}
\numberwithin{equation}{section}
\newtheorem{theorem}{Theorem}[section]
\newtheorem{lemma}[theorem]{Lemma}
\newtheorem{proposition}[theorem]{Proposition}
\newtheorem{corollary}[theorem]{Corollary}
\theoremstyle{definition}
\newtheorem{definition}[theorem]{Definition}
\theoremstyle{remark}

\newcommand{\R}{\mathbb R}
\newcommand{\N}{\mathbb N}
\newcommand{\G}{\mathcal G}
\newcommand{\bdim}{\dim_{\mathrm B}}
\newcommand{\lbdim}{\underline{\dim}_{\mathrm B}}
\newcommand{\ubdim}{\overline{\dim}_{\mathrm B}}
\newcommand{\osc}{\operatorname{osc}}
\newcommand{\Var}{\operatorname{Var}}

\newcommand{\Had}{\mathcal H}
\newcommand{\Weyl}{\mathcal W}
\newcommand{\norm}[1]{\left\lVert#1\right\rVert}

\newcommand{\dd}{\,\mathrm d}
\newcommand{\doi}[1]{\href{https://doi.org/#1}{\nolinkurl{#1}}}
\setlist[enumerate]{label=\textup{(\roman*)},leftmargin=2em,itemsep=3pt}

\makeatletter
\renewcommand{\@maketitle}{%
  \newpage\null\vskip 2em
  \begin{center}
    {\LARGE \@title\par}
    \vskip .35em
    {\large \begin{tabular}[t]{c}\@author\end{tabular}\par}
  \end{center}
  \par
}
\makeatother
\title{Sharp box-dimension bounds for fractional integrals\protect\\ under unequal H\"older and anti-H\"older exponents}
\author{%
  \begin{tabular}{@{}c@{\hspace{3em}}c@{}}
    Long Tian & Yongshun Liang\\[-1pt]
    \small\href{mailto:tianlong19850812@163.com}{\nolinkurl{tianlong19850812@163.com}} &
    \small\href{mailto:yongshun_liang@163.com}{\nolinkurl{yongshun_liang@163.com}}
  \end{tabular}\\[0.25em]
  \small School of Mathematics and Statistics, Nanjing University of Science and Technology,\\[-0.1em]
  \small Nanjing 210094, China.
}
\date{}

\begin{document}
\maketitle
\begin{abstract}
Let $0<h\le H<1$ and $0<\sigma<1$. Graph box dimensions are studied for fractional integrals of functions satisfying an $h$-Holder upper bound on increments and an $H$-anti-H \"lder lower bound on oscillation on a nondegenerate subinterval. For the Riemann--Liouville integral, the following inequalities hold:
\[
\begin{aligned}
\max\{1,2-H-\sigma H/h\}
&\le\lbdim\G(I_{0+}^{\sigma}f,[0,1])\\
&\le\ubdim\G(I_{0+}^{\sigma}f,[0,1])
\le\max\{1,2-h-\sigma\}.
\end{aligned}
\]
The lower bound is obtained by balancing the truncation error of a Marchaud inverse against the oscillation of the input. For each $(h,H)$, a single sparse tent construction attains it for all $\sigma\in(0,1)$. The same bounds and their optimality are established for Hadamard integrals on compact intervals bounded away from zero and for left- and right-sided Liouville--Weyl integrals under explicit weighted tail conditions. Zero-mean periodic Weyl integrals with conditionally convergent tails are treated separately. For Hadamard integrals the bounds are transferred by a bi-Lipschitz logarithmic change of variables. For Weyl integrals the transfer relies on invariance of graph box dimensions under continuous perturbations of bounded variation. Local hypotheses and point-dependent anti-H\"older constants are also treated. The exact linear dimension law is recovered for equal exponents, whereas unequal exponents allow a strictly greater loss of graph dimension.
\end{abstract}

\noindent\textbf{Keywords:} Riemann--Liouville integral; Hadamard integral; Weyl integral; box dimension; anti-H\"older oscillation; bounded variation.

\noindent\textbf{Mathematics Subject Classification (2020):} 26A33, 28A80, 26A16.

\section{Introduction}\label{sec:intro}

Fractional integration improves regularity, but its effect on the geometry of a continuous graph depends on how oscillations are distributed. A Holder estimate bounds the size of an increment in terms of its horizontal scale. An anti-H\"older estimate bounds oscillation from below. When the two exponents agree, the input has a single uniform scale of roughness. When they differ, the oscillation observed on an interval may be concentrated in a much narrower subinterval. This separation of scales affects the extent to which fractional integration smooths the graph. A sharp lower bound reflecting the two scales is established for Riemann--Liouville integrals and then extended to Hadamard and Liouville--Weyl integrals. In all three cases the graph is restricted to a compact nondegenerate interval.

Throughout the paper, the graph of a continuous function $u$ over a compact interval $J$ is
\[
\G(u,J)=\{(x,u(x)):x\in J\}.
\]
Its lower and upper box dimensions are denoted by $\underline{\dim}_{\mathrm B}\mathcal G(u,J)$ and $\overline{\dim}_{\mathrm B}\mathcal G(u,J)$, respectively. When they coincide, their common value is written $\dim_{\mathrm B}\mathcal G(u,J)$. The main results concern box dimension and make no assertion about Hausdorff dimension.

For $0<\sigma<1$, the left-sided Riemann--Liouville integral is
\begin{equation}\label{eq:rl}
(I_{a+}^{\sigma}f)(x)
=\frac{1}{\Gamma(\sigma)}\int_a^x(x-t)^{\sigma-1}f(t)\dd t,
\qquad x\ge a.
\end{equation}
The value at $a$ is zero. For $0<h<1$, the notation $f\in C^h(J)$ means that
\begin{equation}\label{eq:holder}
[f]_{h;J}:=\sup_{x\ne y\in J}\frac{|f(x)-f(y)|}{|x-y|^h}<\infty.
\end{equation}
For a subinterval $K$, set $\osc_K f=\sup_K f-\inf_K f$.

\begin{definition}\label{def:anti}
A continuous function is uniformly $H$-anti-H\"older on a nondegenerate interval $J$ if there exist $c>0$ and $r_0>0$ such that
\begin{equation}\label{eq:anti}
\osc_K f\ge c|K|^H
\end{equation}
for every nondegenerate closed interval $K\subset J$ with $|K|\le r_0$.
\end{definition}

Here \emph{anti-H\"older} refers to a quantitative lower bound on local oscillation. No reverse increment inequality is required for every pair of points. The condition is distinct from the reverse H\"older inequalities for weights. A local version with a fixed exponent and point-dependent constants is also considered. A map is called bi-Lipschitz if both the map and its inverse on its image are Lipschitz. Equivalently, distances are bounded above and below by fixed positive multiples of the original distances.

\begin{theorem}[Unequal-exponent bounds]\label{thm:main}
Let $0<h\le H<1$, let $f\in C^h([0,1])$, and suppose that \eqref{eq:anti} holds on some nondegenerate closed interval $J\subset[0,1]$. For every $0<\sigma<1$, put $g=I_{0+}^{\sigma}f$. Then
\begin{equation}\label{eq:main}
\begin{split}
\max\left\{1,\,2-H-\frac{H}{h}\sigma\right\}
&\le \lbdim\G(g,[0,1])\\
&\le \ubdim\G(g,[0,1])
\le \max\{1,\,2-h-\sigma\}.
\end{split}
\end{equation}
In particular, the lower bound is greater than one when
\begin{equation}\label{eq:threshold}
0<\sigma<\frac{h(1-H)}{H}.
\end{equation}
There is no restriction on $f(0)$.
\end{theorem}

The factor $H/h$ converts between two scales. An interval of length $\delta$ carries oscillation of size $\delta^H$, while the $h$-H\"older upper bound allows this oscillation to occur over a distance comparable to $\delta^{H/h}$. Fractional integration can therefore attenuate it by the factor $\delta^{\sigma H/h}$. Liang and Su \cite[Remark~3.13]{LiangSuGeometric} posed the equal-exponent lower-bound problem for $f(0)=0$.

\begin{theorem}[Optimality for every pair of exponents]\label{thm:sharp}
For every $0<h\le H<1$, there exists a nonnegative function $f\in C^h([0,1])$ with $f(0)=f(1)=0$ such that:
\begin{enumerate}
\item $f$ is uniformly $H$-anti-H\"older on $[0,1]$;
\item $\bdim\G(f,[0,1])=2-H$;
\item for every $0<\sigma<1$,
\begin{equation}\label{eq:sharp}
\bdim\G(I_{0+}^{\sigma}f,[0,1])
=\max\left\{1,\,2-H-\frac{H}{h}\sigma\right\};
\end{equation}
\item $f\notin C^\kappa([0,1])$ for every $h<\kappa\le1$.
\end{enumerate}
Moreover, if $H+\sigma H/h>1$, then $I_{0+}^{\sigma}f$ has bounded variation.
\end{theorem}

Consequently the lower bound in Theorem~\ref{thm:main} is the best universal lower bound under its hypotheses. The theorem concerns optimality over a function class, not a claim that all functions in that class have the displayed dimension.

\begin{corollary}[Equal exponents]\label{cor:equal}
If the hypotheses of Theorem~\ref{thm:main} hold with $h=H$, then
\[
\bdim\G(I_{0+}^{\sigma}f,[0,1])
=\max\{1,2-H-\sigma\},\qquad 0<\sigma<1.
\]
\end{corollary}

The following theorem extends the estimates to Hadamard and Weyl integrals, which are defined in Sections~\ref{sec:hadamard} and \ref{sec:weyl}. The signs in the Weyl notation indicate the direction of integration.

\begin{theorem}[Hadamard and Weyl dimension bounds]\label{thm:extensions}
Let $0<h\le H<1$, $0<\sigma<1$, and let $J=[a,b]$ be nondegenerate. Assume $f|_J\in C^h(J)$ and that \eqref{eq:anti} holds on some nondegenerate subinterval of $J$. Let $g$ be any of the following:
\begin{enumerate}
\item $g=\Had_{a+}^{\sigma}f$, the left-sided Hadamard integral, where $0<a<b<\infty$;
\item $g=\Weyl_+^{\sigma}f$, the integral from $-\infty$ to $x$, where $f\in C(( -\infty,b])$ and
\begin{equation}\label{eq:tailintroleft}
\int_{-\infty}^{a}(a-t)^{\sigma-1}|f(t)|\dd t<\infty;
\end{equation}
\item $g=\Weyl_-^{\sigma}f$, the integral from $x$ to $+\infty$, where $f\in C([a,\infty))$ and
\begin{equation}\label{eq:tailintroright}
\int_b^{\infty}(t-b)^{\sigma-1}|f(t)|\dd t<\infty.
\end{equation}
\end{enumerate}
Then $g\in C(J)$ and
\begin{equation}\label{eq:extensionbounds}
\begin{split}
\max\{1,2-H-\sigma H/h\}
&\le\lbdim\G(g,J)\\
&\le\ubdim\G(g,J)\le\max\{1,2-h-\sigma\}.
\end{split}
\end{equation}
In each case the lower bound is optimal for every $(h,H,\sigma)$. When $h=H$, the graph box dimension exists and equals $\max\{1,2-H-\sigma\}$.
\end{theorem}

For Hadamard integrals, optimality holds on every fixed $[a,b]\subset(0,\infty)$. For Weyl integrals, the lower bound is attained even for continuous inputs of compact support, with the H\"older and anti-H\"older conditions imposed on the observation interval. For zero-mean periodic inputs, \eqref{eq:extensionbounds} is established in Proposition~\ref{prop:periodicweyl} without the absolute tail conditions above. Optimality in this periodic class is proved in Corollary~\ref{cor:periodicsharp}.

For Liouville integration, Liu and Liang \cite{LiuLiangWavelet} used wavelets to establish a linear dimension law under their H\"older and anti-H\"older hypotheses. For Riemann--Liouville integration, they proved that the dimension of a regular monofractal of box dimension $s$ drops by $\sigma$ when $0<\sigma<s-1$ \cite{LiuLiangMono}. The classical H\"older regularity estimates go back to Hardy and Littlewood \cite{HL}. For continuous inputs, graph dimensions were studied by Liang and Su \cite{LiangSu2016} and, for graphs of dimension one, by Liang \cite{LiangFractional}. H\"older-based dimension estimates were developed in \cite{LiangHolder,WuUpper}. Wu \cite{WuLinearity} established a linear law under additional hypotheses. Wu, Ji, and Zhang \cite{WuJiZhang} obtained an upper bound in terms of the upper box dimension of the input graph.

For Hadamard integration, Tian \cite{TianHadamard} estimated the H\"older index and box dimension, and Xiao \cite{XiaoHadamard} compared the dimensions before and after integration. The logarithmic conjugacy used in Section~\ref{sec:hadamard} is described in \cite{LiuLiuHadamard}. For Weyl integrals, H\"older regularity and box dimension were studied by Tian \cite{TianWeyl}, upper box dimension as a function of the integration order by Gao, Liang, and Xiao \cite{GaoWeyl}, and inputs of bounded variation by Mu and coauthors \cite{MuWeyl}.

Here the lower bound depends on both the exponent $h$ controlling increments from above and the exponent $H$ controlling oscillation from below. It is attained by the same construction for every integration order. Both box dimensions are preserved by the logarithmic change of variables on compact intervals away from zero. For Weyl integrals, the history outside the observation interval contributes a continuous term of bounded variation. Its removal permits the local Riemann--Liouville estimate to be applied. Periodic inputs are treated separately by cancellation. For dimension estimates involving Weyl--Marchaud derivatives, see \cite{LiangLiuMarchaud}. The truncated Marchaud inverse estimate used here is proved directly.

In Section~\ref{sec:osc} the dimension bounds and their local forms are deduced from covering, smoothing, and inverse estimates. Optimality is proved by the construction in Section~\ref{sec:sharp}. The bounds are extended to Hadamard and Weyl integrals in Section~\ref{sec:extensions}. Their scope and their relation to results for other function classes and operators are discussed in Section~\ref{sec:scope}.

\section{Box dimension estiamtes}\label{sec:osc}

For $m\in\N_+$, write $\delta=1/m$ and
\[
Q_{k,\delta}=[k\delta,(k+1)\delta],\qquad 0\le k<m.
\]
For continuous $u:[0,1]\to\R$, define
\begin{equation}\label{eq:oscsum}
S_\delta(u)=\sum_{k=0}^{m-1}\osc_{Q_{k,\delta}}u.
\end{equation}
Let $N_\delta(E)$ denote the least number of closed axis-parallel squares of side $\delta$ covering a compact set $E\subset\R^2$. The lower and upper limits of $\log N_\delta(E)/\log(1/\delta)$ define its box dimensions. Restricting $\delta$ to $1/m$ does not change these limits, by monotonicity and the fact that successive scales have bounded ratio tending to one.

\begin{lemma}[Covering by oscillations]\label{lem:cover}
There are two absolute constants $c,C>0$ such that
\begin{equation}\label{eq:cover}
c\delta^{-1}\bigl(1+S_\delta(u)\bigr)
\le N_\delta(\G(u,[0,1]))
\le C\delta^{-1}\bigl(1+S_\delta(u)\bigr).
\end{equation}
For the covering convention above, one may take $c=1/4$ and $C=1$.
\end{lemma}
\begin{proof}
Write $r_k=\osc_{Q_{k,\delta}}u$ and $a_k=\min_{Q_{k,\delta}}u$. Above the horizontal interval $Q_{k,\delta}$, stack closed squares of side $\delta$ starting at height $a_k$. The number
\[
L_k=\max\{1,\lceil r_k/\delta\rceil\}
\le1+r_k/\delta
\]
suffices to cover the graph over that interval, including its endpoints, even when $r_k=0$. Summation over the $m=\delta^{-1}$ columns then yields
\[
N_\delta(\G(u,[0,1]))\le\sum_{k=0}^{m-1}L_k
\le\delta^{-1}\bigl(1+S_\delta(u)\bigr).
\]

Conversely, take any cover of the graph by $N$ closed squares of side $\delta$, with arbitrary positions. Let $n_k$ count the covering squares whose horizontal projections intersect $Q_{k,\delta}$. Continuity makes the vertical projection of the graph over $Q_{k,\delta}$ the full interval $[a_k,a_k+r_k]$. The vertical projections of those $n_k$ squares cover this interval, so $r_k\le\delta n_k$. A closed interval of length $\delta$ intersects at most three of the closed mesh intervals: the possible third intersection occurs only at a grid endpoint. Hence
\[
S_\delta(u)\le\delta\sum_{k=0}^{m-1}n_k\le3\delta N.
\]
Since the horizontal projections of the same squares cover $[0,1]$, one also has $1\le\delta N$. Adding the two inequalities yields $1+S_\delta(u)\le4\delta N$. The lower bound follows on minimizing over covers, with intersections at mesh boundaries included in the count. For the upper bound, squares may be freely positioned by definition of $N_\delta$. Restriction to a fixed lattice changes only the covering constants, so both box dimensions are preserved.
\end{proof}

\begin{lemma}[Bounded-variation invariance]\label{lem:bv}
If $u,v$ are continuous on a compact nondegenerate interval and $v$ has bounded variation, then adding $v$ to $u$ leaves both graph box dimensions of $u$ and $v$ unchanged.
\end{lemma}
\begin{proof}
It suffices to work on $[0,1]$, because affine changes in the horizontal coordinate are bi-Lipschitz. For every mesh,
\[
S_\delta(u+v)\le S_\delta(u)+S_\delta(v)
\le S_\delta(u)+\Var(v).
\]
Applying the same argument to $u=(u+v)-v$ shows that
\[
\bigl|S_\delta(u+v)-S_\delta(u)\bigr|\le\Var(v).
\]
Thus $1+S_\delta(u+v)$ and $1+S_\delta(u)$ are comparable by constants independent of $\delta$. Lemma~\ref{lem:cover} yields the assertion for both limits. In particular, a continuous bounded-variation function has graph box dimension one.
\end{proof}

\begin{lemma}[Lower oscillation mass]\label{lem:mass}
If \eqref{eq:anti} holds on a nondegenerate interval $J$, then, for all sufficiently small mesh sizes,
\begin{equation}\label{eq:mass}
S_\delta(f)\ge c_J\delta^{H-1}
\end{equation}
with a constant $c_J>0$.
\end{lemma}
\begin{proof}
For $\delta\le |J|/4$, at least $|J|/(2\delta)$ complete mesh intervals lie in $J$. If also $\delta\le r_0$, each contributes at least $c\delta^H$. Their sum proves the claim.
\end{proof}

\begin{lemma}[Translated meshes]\label{lem:translate}
Suppose $u\in C([0,1])$, $u(0)=0$, and extend $u$ by zero to $(-\infty,0)$. Then, for every $t\ge0$,
\begin{equation}\label{eq:translate}
\sum_{k=0}^{m-1}\osc_{Q_{k,\delta}-t}u\le 2S_\delta(u).
\end{equation}
\end{lemma}
\begin{proof}
An interval of the translated mesh meets at most two original mesh intervals in their interiors. By continuity, its oscillation is bounded by the sum of their oscillations. Each original mesh interval is used at most twice. Intervals wholly in the negative half-line contribute zero. An interval crossing zero is handled in the same way because the zero extension is continuous at zero. None of the translated intervals extends to the right of one.
\end{proof}

The lower bound requires the following analytic estimates with the stated endpoint normalization. For $\alpha,\beta>0$ and a continuous function $F$, the semigroup identity
\begin{equation}\label{eq:semigroup}
I_{0+}^{\alpha}I_{0+}^{\beta}F=I_{0+}^{\alpha+\beta}F
\end{equation}
follows from Fubini's theorem and
\[
\int_t^x(x-s)^{\alpha-1}(s-t)^{\beta-1}\dd s
=(x-t)^{\alpha+\beta-1}B(\alpha,\beta).
\]
Absolute integrability follows from boundedness of $F$ and positivity and local integrability of the power kernels.

\begin{lemma}[Normalized smoothing]\label{lem:smoothing}
Let $0<h<1$, $0<\sigma<1$, and $\beta=h+\sigma<1$. If $F\in C^h([0,1])$ and $F(0)=0$, then $G=I_{0+}^{\sigma}F$, extended by zero to the negative half-line, is $\beta$-H\"older on $(-\infty,1]$. More precisely,
\begin{equation}\label{eq:smoothing}
[G]_{\beta;(-\infty,1]}\le C_{h,\sigma}[F]_{h;[0,1]}.
\end{equation}
\end{lemma}
Here $[G]_{\beta;(-\infty,1]}$ means the beta-th order semi-norm of $G$ on $(-\infty,1]$, and $[F]_{h;[0,1]}$ means the $h$-th order semi-norm of $F$ on $[0,1]$.
\begin{proof}
Set $M=[F]_h$. Since $F(0)=0$, one has $|F(t)|\le Mt^h$. The beta integral then yields
\begin{equation}\label{eq:initial-growth}
|G(x)|\le M\frac{\Gamma(h+1)}{\Gamma(h+\sigma+1)}x^{h+\sigma}.
\end{equation}
Take $0\le x<x+r\le1$. Adding and subtracting $F(x)$ in the two integrals shows
\begin{align*}
\Gamma(\sigma)\bigl(G(x+r)-G(x)\bigr)
={}&\int_0^x\bigl[(x+r-t)^{\sigma-1}-(x-t)^{\sigma-1}\bigr]
       (F(t)-F(x))\dd t\\
&+\int_x^{x+r}(x+r-t)^{\sigma-1}(F(t)-F(x))\dd t\\
&+\frac{F(x)}{\sigma}\bigl((x+r)^\sigma-x^\sigma\bigr).
\end{align*}
The absolute value of the first integral is at most
\[
M\int_0^x\bigl[s^{\sigma-1}-(s+r)^{\sigma-1}\bigr]s^h\dd s
\le Mr^{h+\sigma}\int_0^\infty
\bigl[z^{\sigma-1}-(z+1)^{\sigma-1}\bigr]z^h\dd z.
\]
The last integral is finite: its integrand is integrable near zero since $h+\sigma>0$, and is $O(z^{h+\sigma-2})$ at infinity since $h+\sigma<1$. The second integral is bounded by
\[
M\int_0^r(r-s)^{\sigma-1}s^h\dd s=C_{h,\sigma}Mr^{h+\sigma}.
\]
For the last term, if $x\le r$, use $|F(x)|\le Mr^h$ and $(x+r)^\sigma-x^\sigma\le r^\sigma$. When $x\ge r$, the following estimate is obtained by the mean value theorem:
\[
|F(x)|\bigl((x+r)^\sigma-x^\sigma\bigr)
\le M\sigma r x^{h+\sigma-1}\le M\sigma r^{h+\sigma}.
\]
This proves the estimate when both points are nonnegative. If one point is negative, \eqref{eq:initial-growth} proves the required estimate across zero. If both are negative, the increment vanishes.
\end{proof}

The Marchaud derivative is a fractional derivative represented by a weighted integral of backward increments, rather than by a classical derivative of the input. In the normalization below it recovers $F$ from its fractional integral $G$.

\begin{lemma}[Marchaud inversion]\label{lem:inverse}
Under the assumptions of Lemma~\ref{lem:smoothing}, for $0\le x\le1$,
\begin{equation}\label{eq:marchaud}
F(x)=\frac{\sigma}{\Gamma(1-\sigma)}
\int_0^\infty\frac{G(x)-G(x-t)}{t^{1+\sigma}}\dd t.
\end{equation}
The integral is absolutely convergent.
\end{lemma}
\begin{proof}
Near zero, the integrand is bounded by a constant times $t^{h-1}$ by Lemma~\ref{lem:smoothing}. For $t\ge1$, boundedness and the zero extension imply an integrable bound of order $t^{-1-\sigma}$.

First suppose $F$ is smooth and vanishes on a neighborhood to the left and right of zero. Then $G$ is smooth across zero. Integration by parts, with vanishing boundary terms at zero and infinity, shows
\[
\frac{\sigma}{\Gamma(1-\sigma)}
\int_0^\infty\frac{G(x)-G(x-t)}{t^{1+\sigma}}\dd t
=\frac{1}{\Gamma(1-\sigma)}\int_0^\infty t^{-\sigma}G'(x-t)\dd t
=I_{0+}^{1-\sigma}G'(x).
\]
Since $F(0)=0$, differentiation of the convolution yields $G'=I_{0+}^{\sigma}F'$. By \eqref{eq:semigroup}, the last expression is $I_{0+}^1F'=F$.

For a general $F$, extend it by zero to the negative half-line and use a smooth nonnegative kernel $\rho_\eta$ of integral one supported in $(\eta/2,\eta)$. Define the causal mollifications
\[
F_\eta(x):=\int_0^\infty\rho_\eta(s)F(x-s)\dd s,
\qquad x\le1.
\]
They are smooth on the relevant interval, vanish near zero, and converge uniformly to $F$ on $[0,1]$. Their $h$-H\"older seminorms on $(-\infty,1]$ are uniformly bounded. A harmless extension to the right can be used to define the convolution at one; it does not enter any left-sided integral below. By Fubini's theorem,
\[
G_\eta:=I_{0+}^{\sigma}F_\eta=\rho_\eta*G
\]
on $(-\infty,1]$. Thus $G_\eta\to G$ uniformly there, and their $\beta$-H\"older seminorms are uniformly bounded. Apply the smooth identity and pass to the limit. For $0<t<1$ the integrands are dominated by $Ct^{h-1}$, and for $t\ge1$ by $Ct^{-1-\sigma}$. Dominated convergence proves \eqref{eq:marchaud}; the identity at zero is immediate.
\end{proof}

For $\varepsilon>0$, define the truncated inverse
\begin{equation}\label{eq:truncated}
T_\varepsilon G(x)=\frac{\sigma}{\Gamma(1-\sigma)}
\int_\varepsilon^\infty
\frac{G(x)-G(x-t)}{t^{1+\sigma}}\dd t.
\end{equation}

\begin{lemma}[Truncation error and oscillation bound]\label{lem:trunc}
Under the assumptions of Lemma~\ref{lem:smoothing}, there are constants $A,B$ independent of $\varepsilon\in(0,1]$ and $\delta=1/m$ such that
\begin{align}
\norm{F-T_\varepsilon G}_{\infty;[0,1]}&\le A\varepsilon^h,
\label{eq:error}\\
S_\delta(T_\varepsilon G)&\le B\varepsilon^{-\sigma}S_\delta(G).
\label{eq:Tosc}
\end{align}
Consequently,
\begin{equation}\label{eq:key}
S_\delta(F)\le B\varepsilon^{-\sigma}S_\delta(G)
+2A\delta^{-1}\varepsilon^h.
\end{equation}
\end{lemma}
\begin{proof}
Subtract \eqref{eq:truncated} from \eqref{eq:marchaud}. The omitted part is bounded by
\[
\frac{\sigma[G]_\beta}{\Gamma(1-\sigma)}
\int_0^\varepsilon t^{\beta-1-\sigma}\dd t
=\frac{\sigma[G]_\beta}{h\Gamma(1-\sigma)}\varepsilon^h,
\]
which proves \eqref{eq:error}. For any mesh interval $Q$, it follows from \eqref{eq:truncated}, by taking differences, that
\[
\osc_Q(T_\varepsilon G)
\le\frac{\sigma}{\Gamma(1-\sigma)}
\int_\varepsilon^\infty t^{-1-\sigma}
\bigl(\osc_QG+\osc_{Q-t}G\bigr)\dd t.
\]
On summing over the mesh intervals and applying Lemma~\ref{lem:translate}, estimate \eqref{eq:Tosc} is obtained with $B=3/\Gamma(1-\sigma)$. On each of the $m$ intervals the oscillation of the error is at most twice its uniform norm. Estimate \eqref{eq:key} follows.
\end{proof}

To prove Theorem~\ref{thm:main}, the inverse estimate is combined with the lower oscillation bound. The contribution of the endpoint value is treated separately.

\begin{proof}[Proof of Theorem~\ref{thm:main}]
Set
\[
F=f-f(0),\qquad G=I_{0+}^{\sigma}F.
\]
Then $F(0)=0$, its H\"older seminorm equals that of $f$, and all oscillations are unchanged. Moreover,
\begin{equation}\label{eq:initial-bv}
g(x)=G(x)+\frac{f(0)}{\Gamma(1+\sigma)}x^\sigma.
\end{equation}
The second term is continuous and of bounded variation. By Lemma~\ref{lem:bv}, $g$ and $G$ have the same lower and upper graph box dimensions.

Every continuous graph over $[0,1]$ has lower box dimension at least one. Hence the asserted lower bound is immediate if
\[
E:=H+\frac{H}{h}\sigma\ge1.
\]
Suppose $E<1$. Then $h+\sigma=(h/H)E<1$, so Lemmas~\ref{lem:smoothing}--\ref{lem:trunc} apply. Combining Lemma~\ref{lem:mass} with \eqref{eq:key} yields
\begin{equation}\label{eq:balance}
c_J\delta^{H-1}\le B\varepsilon^{-\sigma}S_\delta(G)
+2A\delta^{-1}\varepsilon^h.
\end{equation}
Choose a fixed $0<\eta\le1$ such that $2A\eta^h\le c_J/2$, and put
\begin{equation}\label{eq:epsilonchoice}
\varepsilon=\eta\delta^{H/h}.
\end{equation}
The last term in \eqref{eq:balance} is then at most $c_J\delta^{H-1}/2$. Thus
\begin{equation}\label{eq:Gmass}
S_\delta(G)\ge\frac{c_J\eta^\sigma}{2B}
\delta^{H-1+\sigma H/h}=c_*\delta^{E-1}.
\end{equation}
By Lemma~\ref{lem:cover}, $N_\delta(\mathcal G(G,[0,1]))\ge c\delta^{E-2}$ for every sufficiently small mesh size. The lower bound is obtained by taking the lower limit.

For the upper bound, if $h+\sigma<1$, Lemma~\ref{lem:smoothing} implies $S_\delta(G)\le C\delta^{h+\sigma-1}$, and hence $\ubdim\G(G,[0,1])\le2-h-\sigma$. If $h+\sigma\ge1$, take any $0<h'<1-\sigma$ with $h'<h$. Since $F\in C^{h'}$, the upper bound $2-h'-\sigma$ follows as before. Letting $h'\uparrow1-\sigma$ reduces this bound to one. Together with \eqref{eq:initial-bv}, this proves \eqref{eq:main} in all cases.
\end{proof}

\begin{proof}[Proof of Corollary~\ref{cor:equal}]
When $h=H$, the lower and upper endpoints in \eqref{eq:main} coincide. The formula also covers $H+\sigma=1$ by the same limiting argument for the upper bound.
\end{proof}

In the preceding proof, the anti-H\"older hypothesis is used only to establish \eqref{eq:mass}. The following more flexible statement is therefore obtained.

\begin{proposition}[Aggregate oscillation hypothesis]\label{prop:aggregate}
The conclusion of Theorem~\ref{thm:main} remains valid if \eqref{eq:anti} is replaced by the assumption that $S_\delta(f)\ge c\delta^{H-1}$ for every sufficiently small $\delta=1/m$.
\end{proposition}
\begin{proof}
This assumption is exactly the lower estimate used in \eqref{eq:balance}; all subsequent steps are unchanged.
\end{proof}

\begin{proposition}[Localization of both hypotheses]\label{prop:local}
Let $f\in C([0,1])$, let $J=[a,b]\subset[0,1]$ be nondegenerate, and assume $f|_J\in C^h(J)$. Suppose \eqref{eq:anti} holds on a nondegenerate subinterval of $J$. Then
\begin{equation}\label{eq:local}
\begin{split}
\max\{1,2-H-\sigma H/h\}
&\le\lbdim\G(I_{0+}^{\sigma}f,J)\\
&\le\ubdim\G(I_{0+}^{\sigma}f,J)
\le\max\{1,2-h-\sigma\}.
\end{split}
\end{equation}
In particular, the same lower bound holds for the graph over $[0,1]$.
\end{proposition}
\begin{proof}
For $x\in[a,b]$, split
\[
I_{0+}^{\sigma}f(x)=I_{a+}^{\sigma}(f|_J)(x)+R_a(x),\qquad
R_a(x)=\frac1{\Gamma(\sigma)}\int_0^a(x-t)^{\sigma-1}f(t)\dd t.
\]
If $a=0$, the history term is zero. For $a>0$, it is continuous on $[a,b]$. This follows from dominated convergence, since for $x\ge a$ the absolute value of the integrand is bounded by $\|f\|_\infty(a-t)^{\sigma-1}$. For $a=x_0<\cdots<x_N=b$, monotonicity of $x\mapsto(x-t)^{\sigma-1}$ implies
\begin{align*}
\sum_{j=1}^N|R_a(x_j)-R_a(x_{j-1})|
&\le\frac1{\Gamma(\sigma)}\int_0^a|f(t)|
\bigl((a-t)^{\sigma-1}-(b-t)^{\sigma-1}\bigr)\dd t\\
&\le\frac{\norm f_\infty a^\sigma}{\Gamma(1+\sigma)}.
\end{align*}
Thus $R_a$ has bounded variation. Lemma~\ref{lem:bv} removes it without changing either box dimension. Rescale $J$ affinely to $[0,1]$ and apply Theorem~\ref{thm:main}. In detail, for $L=b-a$ and $\widetilde f(y)=f(a+Ly)$,
\[
I_{a+}^{\sigma}f(a+Ly)=L^\sigma I_{0+}^{\sigma}\widetilde f(y).
\]
The changes in both coordinates are bi-Lipschitz, proving \eqref{eq:local}. The final assertion follows from inclusion of the local graph in the full graph.
\end{proof}

\begin{proposition}[Point-dependent anti-H\"older constants]\label{prop:baire}
Let $f\in C^h([0,1])$, and fix $0<h\le H<1$ and a nondegenerate closed interval $J\subset(0,1)$. Suppose that, for every $x\in J$, there exist $c_x,r_x>0$ such that
\begin{equation}\label{eq:pointwiseanti}
\osc_{[x-r,x+r]}f\ge c_x r^H
\end{equation}
whenever $0<r<r_x$ and the displayed interval lies in $[0,1]$. Then the lower bound in Theorem~\ref{thm:main} holds. The conclusion remains true if the exponent in \eqref{eq:pointwiseanti} is $H_x\le H$.
\end{proposition}
\begin{proof}
Choose an integer $m_0$ such that $1/m_0<\operatorname{dist}(J,\{0,1\})$. For $m\ge m_0$, let
\[
E_m=\left\{x\in J:
\osc_{[x-r,x+r]}f\ge m^{-1}r^H
\text{ for every }0<r\le m^{-1}\right\}.
\]
For fixed $r$, the oscillation varies continuously under translation of the interval by uniform continuity of $f$. Hence each $E_m$ is closed in $J$. By hypothesis $J=\bigcup_{m\ge m_0}E_m$. The Baire category theorem implies that some $E_m$ contains a nondegenerate subinterval $U$. If $K\subset U$ and $|K|\le2/m$, take $x$ to be its midpoint and $r=|K|/2$. Then
\[
\osc_K f\ge m^{-1}2^{-H}|K|^H.
\]
Theorem~\ref{thm:main} applies. If $H_x\le H$, first observe that $r^{H_x}\ge r^H$ for $0<r<1$ and reduce to the fixed-exponent case.
\end{proof}

\clearpage
\section{A construction attaining the lower bound}\label{sec:sharp}

Now prove Theorem~\ref{thm:sharp}. All estimates are quantitative in the fixed parameters. The construction separates the spacing of oscillations from their width.

Fix $0<h\le H<1$, set
\begin{equation}\label{eq:q}
q=H/h\ge1,
\end{equation}
and choose an even integer $b\ge4$ so large that $b^H\ge4$. Let
\begin{equation}\label{eq:parameters}
d_n=b^{-n},\quad a_n=b^{-Hn},\quad
m_n=\lceil qn\rceil+2,\quad w_n=b^{-m_n},\qquad n\ge1.
\end{equation}
In particular,
\begin{equation}\label{eq:width}
b^{-3}b^{-qn}\le w_n\le b^{-2}b^{-qn},\qquad
w_n\le b^{-2}d_n.
\end{equation}
The rounding of $qn$ is important: $q$ can be any real number at least one, while each width still lies on an integer $b$-adic scale.

Put $\phi(t)=(1-|t|)_+$ and define
\begin{equation}\label{eq:tentseries}
f_n(x):=a_n\sum_{k=0}^{b^n-1}
\phi\left(\frac{x-(k+1/2)d_n}{w_n}\right),
\qquad f(x):=\sum_{n=1}^\infty f_n(x),\quad 0\le x\le1.
\end{equation}
At each level, the supports are disjoint and lie inside their grid cells. Hence $0\le f_n\le a_n$, and $\sum a_n<\infty$ proves uniform convergence, continuity, nonnegativity, and $f(0)=f(1)=0$.

\begin{lemma}[Regularity and uniform oscillation of the input]\label{lem:tentrough}
The function \eqref{eq:tentseries} belongs to $C^h([0,1])$, is uniformly $H$-anti-H\"older on $[0,1]$, and does not belong to $C^\kappa([0,1])$ for any $h<\kappa\le1$.
\end{lemma}
\begin{proof}
Since the Lipschitz constant of $f_n$ is at most $a_n/w_n$,
\begin{equation}\label{eq:tentincrement}
|f_n(x)-f_n(y)|\le a_n\min\{1,|x-y|/w_n\}.
\end{equation}
The case $x=y$ is immediate. For $r=|x-y|>0$, choose the unique integer $N\ge0$ such that $b^{-q(N+1)}<r\le b^{-qN}$. By \eqref{eq:width} and \eqref{eq:tentincrement},
\[
|f(x)-f(y)|\le b^3r\sum_{n=1}^{N}b^{(q-H)n}
+\sum_{n=N+1}^{\infty}b^{-Hn}.
\]
Set $\theta=q-H>0$. The two geometric sums are explicitly
\begin{align*}
\sum_{n=1}^{N}b^{\theta n}
&=\frac{b^\theta(b^{\theta N}-1)}{b^\theta-1}
\le\frac{b^\theta}{b^\theta-1}b^{\theta N},\\
\sum_{n=N+1}^{\infty}b^{-Hn}
&=\frac{b^{-H(N+1)}}{1-b^{-H}}.
\end{align*}
The first identity remains valid for $N=0$, when the sum is empty. The scale inequalities imply $b^{qN}\le r^{-1}$ and $b^{-H(N+1)}<r^{H/q}$. Consequently,
\begin{align*}
r b^{\theta N}&\le r^{1-\theta/q}=r^{H/q},\\
|f(x)-f(y)|&\le
\left(\frac{b^{3+\theta}}{b^\theta-1}
+\frac1{1-b^{-H}}\right)r^{H/q}.
\end{align*}
Since $H/q=h$, this proves the required H\"older estimate at every scale.

For the lower oscillation estimate, fix a level-$n$ tent with center $c=(k+1/2)d_n$ and support $[c-w_n,c+w_n]$. Every coarser function $f_j$, $j<n$, is affine on this support. To verify this, note that its centers belong to $d_n\mathbb Z$, since $b$ is even. Its remaining knots are those centers plus or minus $b^{-m_j}$. If $m_j\le n$, these knots are also in $d_n\mathbb Z$. If $m_j>n$, they lie within distance $d_n/b$ of $d_n\mathbb Z$. In contrast, the level-$n$ support is centered halfway between consecutive points of $d_n\mathbb Z$ and has half-width at most $d_n/b^2$. Since $b\ge4$, it meets none of the coarser knots.

Apply the symmetric second difference at $c$ with step $w_n$. The coarser terms cancel, and the level-$n$ term contributes $a_n$. Each finer term is between zero and its amplitude. Therefore
\begin{equation}\label{eq:tentseconddiff}
f(c)-\frac{f(c-w_n)+f(c+w_n)}2
\ge a_n-\sum_{j>n}a_j
=a_n\left(1-\frac1{b^H-1}\right)
\ge\frac23a_n.
\end{equation}
The oscillation on the support is consequently at least $2a_n/3$.

For any nondegenerate interval $K\subset[0,1]$ of length $r$, choose the least $n\ge1$ for which $2d_n\le r$. Then $K$ contains a complete level-$n$ grid cell and
\[
\frac{r}{2b}<d_n\le\frac r2.
\]
The tent support lies in that cell. Hence, by \eqref{eq:tentseconddiff},
\begin{equation}\label{eq:tentanti}
\osc_K f\ge\frac23a_n
\ge\frac23(2b)^{-H}r^H.
\end{equation}
This is a uniform anti-H\"older estimate.

Finally, if $f\in C^\kappa$ with $\kappa>h$, its oscillation on each of the above supports would be at most $C(2w_n)^\kappa$. By \eqref{eq:width} and \eqref{eq:tentseconddiff}, this would imply that $b^{(q\kappa-H)n}$ is bounded as $n\to\infty$. Since $q\kappa>qh=H$, this is impossible.
\end{proof}

\begin{lemma}[Dimension of the input]\label{lem:tentdim}
The function \eqref{eq:tentseries} satisfies $\bdim\G(f,[0,1])=2-H$.
\end{lemma}
\begin{proof}
The lower bound follows from \eqref{eq:tentanti} and Lemmas~\ref{lem:cover} and \ref{lem:mass}. For the upper bound, $\Var(f_n)=2b^n a_n$ and $\norm{f_n}_\infty=a_n$, so
\[
S_\delta(f_n)\le\min\{2b^na_n,\delta^{-1}a_n\}.
\]
Oscillation is subadditive, and the defining series converges uniformly. Choose $N\ge0$ with $b^{-(N+1)}<\delta\le b^{-N}$. Summing the finite and infinite geometric series yields
\begin{align*}
\sum_{n=1}^{N}b^{(1-H)n}
&=\frac{b^{1-H}\bigl(b^{(1-H)N}-1\bigr)}{b^{1-H}-1}
\le\frac{b^{1-H}}{b^{1-H}-1}\delta^{H-1},\\
\delta^{-1}\sum_{n=N+1}^{\infty}b^{-Hn}
&=\frac{\delta^{-1}b^{-H(N+1)}}{1-b^{-H}}
\le\frac{\delta^{H-1}}{1-b^{-H}}.
\end{align*}
Here $b^N\le\delta^{-1}$ and $b^{-(N+1)}<\delta$ were used separately. Thus
\[
S_\delta(f)\le2\sum_{n=1}^{N}b^{(1-H)n}
+\delta^{-1}\sum_{n=N+1}^{\infty}b^{-Hn}
\le C_{b,H}\delta^{H-1}.
\]
Lemma~\ref{lem:cover} proves the upper bound.
\end{proof}

The next two lemmas control the fractional integral of one level. Subtracting its mean removes the slowly varying accumulated contribution of many positive tents.

\begin{lemma}[Fractional integral of a narrow periodic tent]\label{lem:periodic}
For $0<\eta\le1/4$, let
\[
p_\eta(s)=\sum_{k\in\mathbb Z}\phi\left(\frac{s-k-1/2}{\eta}\right),
\qquad v_\eta(s)=p_\eta(s)-\eta.
\]
Then $v_\eta$ is one-periodic, has zero mean, and, for $0<\sigma<1$,
\begin{equation}\label{eq:periodic}
\sup_{X\ge0}\left|\frac1{\Gamma(\sigma)}
\int_0^X(X-s)^{\sigma-1}v_\eta(s)\dd s\right|
\le C_\sigma\eta^\sigma.
\end{equation}
\end{lemma}
\begin{proof}
Since $\int_\R\phi=1$, the mean of $p_\eta$ is $\eta$. Moreover,
\[
\int_0^1|v_\eta(s)|\dd s\le2\eta.
\]
Its primitive $V_\eta(s)=\int_0^s v_\eta(t)\dd t$ is one-periodic and satisfies $V_\eta(0)=0$ and $\norm{V_\eta}_\infty\le2\eta$.

Fix $X\ge0$ and put $A=\max\{0,X-2\}$. A tent support $[k+1/2-\eta,k+1/2+\eta]$ can intersect $[A,X]$ only if its center belongs to $[A-\eta,X+\eta]$. This latter interval has length at most $2+2\eta\le5/2$, using the assumption $\eta\le1/4$. Since distinct centers are one unit apart, it contains at most $\lfloor 5/2\rfloor+1=3$ centers. Thus the union of the tent supports inside $[A,X]$ has length at most $6\eta$. The same assumption makes the supports disjoint, so $0\le p_\eta\le1$.

For a measurable set of length $\ell$ in $[A,X]$, monotonicity of the kernel shows that its integral against $(X-s)^{\sigma-1}$ is at most $\ell^\sigma/\sigma$. Indeed, after setting $t=X-s$, moving the set to $[0,\ell]$ can only increase the integral of the decreasing function $t^{\sigma-1}$. The tent contribution is therefore at most $(6\eta)^\sigma/\sigma$. The subtracted mean contributes at most
\[
\eta\int_A^X(X-s)^{\sigma-1}\dd s
=\frac{\eta}{\sigma}(X-A)^\sigma
\le\frac{2^\sigma}{\sigma}\eta.
\]
It follows that
\[
\int_A^X(X-s)^{\sigma-1}|v_\eta(s)|\dd s
\le C_\sigma(\eta^\sigma+\eta)\le C_\sigma'\eta^\sigma.
\]
If $A>0$, the remaining integral is evaluated by parts:
\begin{align*}
\int_0^A(X-s)^{\sigma-1}v_\eta(s)\dd s
={}&2^{\sigma-1}V_\eta(A)\\
&-(1-\sigma)\int_0^A(X-s)^{\sigma-2}V_\eta(s)\dd s.
\end{align*}
Both terms have absolute value bounded by $C_\sigma\eta$, uniformly in $X$, because $\norm{V_\eta}_\infty\le2\eta$ and
\[
\int_0^A(X-s)^{\sigma-2}\dd s
\le\int_2^\infty t^{\sigma-2}\dd t
=\frac{2^{\sigma-1}}{1-\sigma}<\infty.
\]
All constants in this argument depend only on the fixed order $\sigma$, not on $X$ or $\eta$.
If $A=0$, there is no remaining interval. Combining the estimates and dividing by $\Gamma(\sigma)$ proves the lemma.
\end{proof}

\begin{lemma}[Variation of an integrated tent]\label{lem:tentvar}
Let $\psi(x)=a\phi((x-c)/w)$ with $a,w>0$ and support contained in $(0,1)$. Then
\begin{equation}\label{eq:tentvar}
\Var_{[0,1]}(I_{0+}^{\sigma}\psi)\le C_\sigma aw^\sigma,
\qquad 0<\sigma<1.
\end{equation}
\end{lemma}
\begin{proof}
Extend the tent by zero and define its full left-sided integral
\[
\Phi_\sigma(z)=\frac1{\Gamma(\sigma)}
\int_{-\infty}^z(z-t)^{\sigma-1}\phi(t)\dd t.
\]
From the identity $\phi(t)=(t+1)_+-2t_++(t-1)_+$, the following expression is obtained:
\[
\Phi_\sigma(z)=\frac{(z+1)_+^{\sigma+1}-2z_+^{\sigma+1}+(z-1)_+^{\sigma+1}}
{\Gamma(\sigma+2)}.
\]
Its derivative is
\[
\Phi_\sigma'(z)=\frac{(z+1)_+^{\sigma}-2z_+^{\sigma}+(z-1)_+^{\sigma}}
{\Gamma(\sigma+1)}.
\]
This derivative is bounded on compact intervals, vanishes for $z<-1$, and is $O(z^{\sigma-2})$ as $z\to+\infty$, by the second-difference formula for the function $z^\sigma$. Hence $\Phi_\sigma'\in L^1(\R)$.

Because the support of $\psi$ is positive, the change of variables $t=c+ws$ yields
\[
I_{0+}^{\sigma}\psi(x)=aw^\sigma\Phi_\sigma((x-c)/w).
\]
Its total variation on $[0,1]$ is at most $aw^\sigma\norm{\Phi_\sigma'}_{L^1(\R)}$, proving \eqref{eq:tentvar}.
\end{proof}

\begin{proof}[Proof of Theorem~\ref{thm:sharp}]
Use \eqref{eq:tentseries}. The asserted input properties follow from Lemmas~\ref{lem:tentrough} and \ref{lem:tentdim}. Fix any $0<\sigma<1$ and define
\[
M_n:=\int_0^1 f_n(x)\dd x=b^na_nw_n,
\qquad U_n:=I_{0+}^{\sigma}(f_n-M_n).
\]
The mean of the periodic extension of $f_n$ is also $M_n$. With $\eta_n=w_n/d_n\le b^{-2}$, the centered function can be written as
\[
f_n(x)-M_n=a_n v_{\eta_n}(x/d_n),\qquad 0\le x\le1.
\]
After rescaling, Lemma~\ref{lem:periodic} yields
\begin{equation}\label{eq:Unsup}
\norm{U_n}_{\infty;[0,1]}
\le C_\sigma a_n d_n^\sigma\eta_n^\sigma
=C_\sigma a_nw_n^\sigma.
\end{equation}
By Lemma~\ref{lem:tentvar} and the variation of $x^\sigma$,
\begin{align}
\Var(U_n)
&\le C_\sigma b^na_nw_n^\sigma
+\frac{M_n}{\Gamma(1+\sigma)}
\le C_\sigma' b^na_nw_n^\sigma.
\label{eq:Unvar}
\end{align}
The last inequality uses $w_n\le w_n^\sigma$.

The sums of $\norm{U_n}_\infty$ and of $M_n$ converge. Indeed, the former follows from \eqref{eq:Unsup}, and $M_n=a_n(w_n/d_n)\le a_n$. Uniform convergence of the original input series and continuity of $I_{0+}^{\sigma}:C([0,1])\to C([0,1])$ imply
\begin{equation}\label{eq:gseries}
I_{0+}^{\sigma}f(x)
=U(x)+\frac{M}{\Gamma(1+\sigma)}x^\sigma,
\qquad U=\sum_{n=1}^\infty U_n,\quad M=\sum_{n=1}^\infty M_n.
\end{equation}
For completeness, the operator norm here is at most $1/\Gamma(1+\sigma)$, directly from \eqref{eq:rl}. The last term in \eqref{eq:gseries} is of bounded variation and may be removed when computing box dimensions.

Set $E=H+q\sigma$. Equations \eqref{eq:width}, \eqref{eq:Unsup}, and \eqref{eq:Unvar} imply
\begin{equation}\label{eq:Unosc}
S_\delta(U_n)
\le C_\sigma\min\{b^n a_nw_n^\sigma,\delta^{-1}a_nw_n^\sigma\}
\le C\min\{b^{(1-E)n},\delta^{-1}b^{-En}\}.
\end{equation}
By uniform convergence and subadditivity of oscillation, this can be summed over $n$. Choose $N\ge0$ with $b^{-(N+1)}<\delta\le b^{-N}$. To make the threshold calculation explicit, write
\[
A_N=\sum_{n=1}^{N}b^{(1-E)n},\qquad
B_N=\delta^{-1}\sum_{n=N+1}^{\infty}b^{-En}.
\]
For $E\ne1$, the finite geometric sum is
\[
A_N=\frac{b^{1-E}\bigl(b^{(1-E)N}-1\bigr)}{b^{1-E}-1},
\]
whereas $A_N=N$ if $E=1$. For every $E>0$,
\[
B_N=\frac{\delta^{-1}b^{-E(N+1)}}{1-b^{-E}}
\le\frac{\delta^{E-1}}{1-b^{-E}}.
\]
If $E<1$, use $b^N\le\delta^{-1}$ to obtain $A_N\le C_{b,E}\delta^{E-1}$. If $E=1$, use $N\le\log(1/\delta)/\log b$ and $B_N\le(1-b^{-1})^{-1}$. If $E>1$, the full series bounding $A_N$ converges, while $\delta^{E-1}\le1$. Consequently,
\begin{equation}\label{eq:threecases}
S_\delta(U)\le C\left(
\sum_{n\le N}b^{(1-E)n}
+\delta^{-1}\sum_{n>N}b^{-En}\right)
\le
\begin{cases}
C'\delta^{E-1},&E<1,\\
C'(1+\log(1/\delta)),&E=1,\\
C',&E>1.
\end{cases}
\end{equation}
It follows from Lemma~\ref{lem:cover} that
\[
\ubdim\G(U,[0,1])\le\max\{1,2-E\}.
\]
Theorem~\ref{thm:main} supplies the matching lower bound because $q=H/h$. Lemma~\ref{lem:bv} and \eqref{eq:gseries} prove \eqref{eq:sharp}. The function $f$ was chosen independently of $\sigma$, so the identity holds for every order in $(0,1)$ for the same function.

If $E>1$, \eqref{eq:Unvar} implies $\sum_n\operatorname{Var}(U_n)<\infty$. By uniform convergence, $\operatorname{Var}(U)\le\sum_n\operatorname{Var}(U_n)$, as is seen by applying the triangle inequality on a finite partition and passing to the limit in the partial sums. The bounded-variation assertion now follows from \eqref{eq:gseries}.
\end{proof}

\begin{corollary}[No nontrivial bound from a positive H\"older exponent alone]\label{cor:failure}
Fix $0<h<1$ and $0<\sigma<1$. There exists a continuous nonnegative $f\in C^h([0,1])$, uniformly anti-H\"older with some exponent $H<1$, such that $I_{0+}^{\sigma}f$ has bounded variation and graph box dimension one.
\end{corollary}
\begin{proof}
Choose
\[
\max\left\{h,\frac{h}{h+\sigma}\right\}<H<1.
\]
The interval is nonempty. Since $H+\sigma H/h>1$, the example follows from Theorem~\ref{thm:sharp}.
\end{proof}

\section{Extention to Hadamard and Weyl integrals}\label{sec:extensions}
\subsection{Hadamard integrals}\label{sec:hadamard}

For $0<a<b<\infty$ and $0<\sigma<1$, define
\begin{equation}\label{eq:haddefinition}
(\Had_{a+}^{\sigma}f)(x)
:=\frac1{\Gamma(\sigma)}\int_a^x
\left(\log\frac{x}{t}\right)^{\sigma-1}f(t)\frac{\dd t}{t},
\qquad a\le x\le b,
\end{equation}
with value zero at $a$. This is the Hadamard fractional integral, not a Hadamard finite-part regularization. Its logarithmic conjugacy with the Riemann--Liouville integral is classical; see, for example, \cite[Eq.~(7)]{LiuLiuHadamard}. The change of variables and the geometric estimates are recorded here because both are needed for the dimension result.

\begin{lemma}[Logarithmic transfer]\label{lem:logtransfer}
Put $L=\log(b/a)$ and $F(s)=f(ae^s)$ for $0\le s\le L$. Then
\begin{equation}\label{eq:hadconjugacy}
(\Had_{a+}^{\sigma}f)(ae^s)=I_{0+}^{\sigma}F(s).
\end{equation}
The maps $f\mapsto F$ and its inverse preserve the H\"older exponent and the interval anti-H\"older exponent, with possibly different constants. Moreover, the lower and upper box dimensions of the graph on the left over $[a,b]$ equal the respective dimensions of the graph on the right over $[0,L]$.
\end{lemma}
\begin{proof}
Substitute $t=ae^u$ in \eqref{eq:haddefinition}. Then $\dd t/t=\dd u$ and $\log(ae^s/t)=s-u$, proving \eqref{eq:hadconjugacy}. For $s,t\in[0,L]$, the mean value theorem implies
\begin{equation}\label{eq:expbilipschitz}
a|s-t|\le|ae^s-ae^t|\le b|s-t|.
\end{equation}
Thus $[F]_h\le b^h[f]_h$. Conversely, $[f]_h\le a^{-h}[F]_h$.

If $K\subset[0,L]$ is an interval, its image $K_x=ae^K$ is an interval and
\[
a|K|\le|K_x|\le b|K|,\qquad \osc_KF=\osc_{K_x}f.
\]
Suppose $K_x$ lies in the prescribed anti-H\"older interval and $|K|\le r_0/b$. Then $\operatorname{osc}_{K_x}f\ge c|K_x|^H$ implies $\operatorname{osc}_K F\ge ca^H|K|^H$. The converse is obtained by applying the inverse change of variables.

The map $(s,y)\mapsto(ae^s,y)$ is bi-Lipschitz on $[0,L]\times\R$ by \eqref{eq:expbilipschitz}. A bi-Lipschitz map compares covering numbers at scales differing by fixed multiplicative factors, so it preserves both box dimensions. It maps the graph of $I_{0+}^{\sigma}F$ onto that of $\Had_{a+}^{\sigma}f$, proving the last assertion.
\end{proof}

\begin{proof}[Proof of Theorem~\ref{thm:extensions}, Hadamard case]
Lemma~\ref{lem:logtransfer} transfers both input assumptions to $F$ on $[0,L]$. The bounds follow from Theorem~\ref{thm:main} after this interval is rescaled to $[0,1]$. Indeed, for $F_L(u)=F(Lu)$,
\[
I_{0+}^{\sigma}F(Lu)=L^\sigma I_{0+}^{\sigma}F_L(u),
\]
and the constant scalings of the two graph coordinates preserve box dimension. Lemma~\ref{lem:logtransfer} transfers the result back to $[a,b]$. The same identity proves continuity, including at $a$. When $h=H$, the two bounds coincide. Optimality is proved below.
\end{proof}

\begin{proposition}[Hadamard sharpness]\label{prop:hadsharp}
Fix $0<a<b<\infty$ and $0<h\le H<1$. There is a nonnegative $f\in C^h([a,b])$, uniformly $H$-anti-H\"older on $[a,b]$, with $f(a)=f(b)=0$, such that for every $0<\sigma<1$,
\begin{equation}\label{eq:hadsharp}
\bdim\G(\Had_{a+}^{\sigma}f,[a,b])
=\max\{1,2-H-\sigma H/h\}.
\end{equation}
The input has box dimension $2-H$ and does not belong to $C^\kappa([a,b])$ for $h<\kappa\le1$. If $H+\sigma H/h>1$, the integral has bounded variation.
\end{proposition}
\begin{proof}
Let $f_*$ be the function in Theorem~\ref{thm:sharp}. Set $L=\log(b/a)$ and
\[
f(x)=f_*\left(\frac{\log(x/a)}L\right).
\]
Lemma~\ref{lem:logtransfer} and affine scaling transfer all input assertions. For $0\le u\le1$,
\begin{equation}\label{eq:hadexampleidentity}
\Had_{a+}^{\sigma}f(ae^{Lu})
=L^\sigma I_{0+}^{\sigma}f_*(u).
\end{equation}
By Theorem~\ref{thm:sharp} and bi-Lipschitz invariance, \eqref{eq:hadsharp} holds simultaneously for all $\sigma$. Total variation is preserved under increasing reparametrization and is multiplied by $L^\sigma$ under the displayed vertical scaling. The final assertion follows.
\end{proof}

The same estimates admit localization. If $0<a_0\le a<b$, $f\in C([a_0,b])$, and both quantitative hypotheses hold on $[a,b]$ as in Theorem~\ref{thm:extensions}, then \eqref{eq:extensionbounds} holds for $\Had_{a_0+}^{\sigma}f$ restricted to $[a,b]$. To verify this, use $s=\log(x/a_0)$ on the full interval and apply Proposition~\ref{prop:local} after affine rescaling. The history preceding $\log(a/a_0)$ is a continuous bounded-variation term in logarithmic coordinates and remains so after the increasing change back to $x$. Proposition~\ref{prop:baire} also carries over: logarithmic coordinates preserve a uniform anti-H\"older estimate once that proposition has produced a suitable subinterval in the original coordinates.

For completeness, the right-sided operator
\[
\Had_{b-}^{\sigma}f(x)
=\frac1{\Gamma(\sigma)}\int_x^b
\left(\log\frac{t}{x}\right)^{\sigma-1}f(t)\frac{\dd t}{t}
\]
has the same bounds. With $s=\log(b/x)$ and $F(s)=f(be^{-s})$, its value at $be^{-s}$ equals $I_{0+}^{\sigma}F(s)$. This decreasing change of coordinates is bi-Lipschitz on the compact intervals in question, so the preceding proof applies. Replacing $\log(x/a)$ by $\log(b/x)$ in the sharp example proves optimality for the right-sided operator as well.

\subsection{Weyl integrals}\label{sec:weyl}

The two directions of Liouville--Weyl integration are distinguished by the notation
\begin{align}
\Weyl_+^{\sigma}f(x)
&=\frac1{\Gamma(\sigma)}\int_{-\infty}^{x}
(x-t)^{\sigma-1}f(t)\dd t,\label{eq:weylleft}\\
\Weyl_-^{\sigma}f(x)
&=\frac1{\Gamma(\sigma)}\int_x^{\infty}
(t-x)^{\sigma-1}f(t)\dd t.\label{eq:weylright}
\end{align}
The absolute convergence conditions are stated explicitly below. The term Weyl is also used for periodic fractional integrals; that convention is treated separately at the end of this section. No absolute convergence is assumed there.

\begin{lemma}[Absolutely convergent history]\label{lem:weyltail}
Let $a<b$, let $0<\sigma<1$, and suppose a measurable $f$ on $(-\infty,a)$ satisfies \eqref{eq:tailintroleft}. Then
\[
R_a(x)=\frac1{\Gamma(\sigma)}
\int_{-\infty}^{a}(x-t)^{\sigma-1}f(t)\dd t,
\qquad x\in[a,b],
\]
is continuous and has bounded variation. More precisely,
\begin{equation}\label{eq:weyltailvariation}
\Var_{[a,b]}R_a
\le\frac1{\Gamma(\sigma)}\int_{-\infty}^{a}|f(t)|
\bigl[(a-t)^{\sigma-1}-(b-t)^{\sigma-1}\bigr]\dd t<\infty.
\end{equation}
\end{lemma}
\begin{proof}
For $x\ge a$ and $t<a$, the integrand is dominated in absolute value by $(a-t)^{\sigma-1}|f(t)|$, which is integrable. Dominated convergence proves continuity, including at $a$. For a partition $a=x_0<\cdots<x_N=b$, the triangle inequality and monotonicity of the kernel imply
\begin{align*}
\sum_{j=1}^N|R_a(x_j)-R_a(x_{j-1})|
&\le\frac1{\Gamma(\sigma)}\int_{-\infty}^{a}|f(t)|
\sum_{j=1}^N\bigl[(x_{j-1}-t)^{\sigma-1}-(x_j-t)^{\sigma-1}\bigr]\dd t\\
&=\frac1{\Gamma(\sigma)}\int_{-\infty}^{a}|f(t)|
\bigl[(a-t)^{\sigma-1}-(b-t)^{\sigma-1}\bigr]\dd t.
\end{align*}
The last integral is finite by \eqref{eq:tailintroleft}. Taking the supremum over partitions proves \eqref{eq:weyltailvariation}.
\end{proof}

\begin{proof}[Proof of Theorem~\ref{thm:extensions}, Weyl cases]
For $x\in[a,b]$, the left-sided integral converges absolutely and may therefore be split at $a$:
\begin{equation}\label{eq:weylsplit}
\Weyl_+^{\sigma}f(x)=I_{a+}^{\sigma}(f|_{[a,b]})(x)+R_a(x).
\end{equation}
The first term is continuous because its input is continuous, and Lemma~\ref{lem:weyltail} treats the second. By Lemma~\ref{lem:bv}, the two sides of
\begin{equation}\label{eq:weyldimensionidentity}
\dim_*\G(\Weyl_+^{\sigma}f,[a,b])
=\dim_*\G(I_{a+}^{\sigma}f,[a,b])
\end{equation}
are equal for each of $\dim_*=\lbdim$ and $\dim_*=\ubdim$. Apply Theorem~\ref{thm:main} after affine rescaling of $[a,b]$ to obtain \eqref{eq:extensionbounds}.

For the right-sided integral, set $\check f(s)=f(-s)$. Under the substitution $t=-u$,
\[
\Weyl_-^{\sigma}f(x)=\Weyl_+^{\sigma}\check f(-x).
\]
The local hypotheses transfer to $[-b,-a]$, and \eqref{eq:tailintroright} becomes the left-tail condition at $-b$. Reflection preserves both graph box dimensions. The left-sided result therefore proves the right-sided one. When $h=H$, the lower and upper bounds coincide. Proposition~\ref{prop:weylsharp} proves optimality.
\end{proof}

The argument requires H\"older and anti-H\"older control only on the observation interval. Outside it, the tail condition supplies exactly the integrability needed to remove the history. In particular, a change of the input on $(-\infty,a)$ that preserves the weighted tail condition does not change either graph box dimension on $[a,b]$ for the left-sided operator. This is an equality of dimensions, not an equality of the functions.

\begin{proposition}[Weyl sharpness with compact support]\label{prop:weylsharp}
For every $0<h\le H<1$, there exists a nonnegative, compactly supported $f\in C^h(\R)$, uniformly $H$-anti-H\"older on $[0,1]$, such that for every $0<\sigma<1$,
\begin{equation}\label{eq:weylsharp}
\bdim\G(\Weyl_+^{\sigma}f,[0,1])
=\max\{1,2-H-\sigma H/h\}.
\end{equation}
Its input graph on $[0,1]$ has dimension $2-H$, its global H\"older exponent is exactly $h$ in the sense of Theorem~\ref{thm:sharp}, and the displayed integral has bounded variation on $[0,1]$ if $H+\sigma H/h>1$. An analogous example exists for $\Weyl_-^{\sigma}$ and on any prescribed compact observation interval.
\end{proposition}
\begin{proof}
Extend the function $f_*$ of Theorem~\ref{thm:sharp} by zero outside $[0,1]$. Since $f_*(0)=f_*(1)=0$, this extension is globally $h$-H\"older: for a pair with one point outside the interval, compare the inside point with the intervening endpoint. Pairs with both points outside have equal zero values. All input assertions now follow from that theorem, and the tails are absolutely integrable for every $\sigma\in(0,1)$.

For $x\in[0,1]$,
\[
\Weyl_+^{\sigma}f(x)=I_{0+}^{\sigma}f_*(x).
\]
Thus \eqref{eq:weylsharp} and the bounded-variation assertion are proved. The right-sided result follows by reflection, with input $x\mapsto f(1-x)$. Examples on any $[a,b]$ are obtained by translation and positive affine scaling. The resulting nonzero vertical scaling preserves dimension.
\end{proof}

For periodic inputs the weighted absolute tail condition usually fails. The appropriate replacement is cancellation at infinity. The next lemma records precisely what is needed.

\begin{lemma}[A bounded conditionally convergent history]\label{lem:conditionaltail}
Suppose $f$ is measurable and bounded on $(-\infty,a]$, and the improper integral
\[
R_a(a)=\frac1{\Gamma(\sigma)}
\lim_{A\to-\infty}\int_A^a(a-t)^{\sigma-1}f(t)\dd t
\]
exists. Then the same improper history $R_a(x)$ exists for all $x\in[a,b]$, is continuous, and satisfies
\begin{equation}\label{eq:conditionalvariation}
\Var_{[a,b]}R_a
\le\frac{\norm f_{\infty;(-\infty,a]}}{\Gamma(1+\sigma)}(b-a)^\sigma.
\end{equation}
\end{lemma}
\begin{proof}
For $L\ge0$, direct integration followed by passage to the limit yields
\begin{equation}\label{eq:kerneldifference}
\int_0^\infty\bigl[s^{\sigma-1}-(s+L)^{\sigma-1}\bigr]\dd s
=\frac{L^\sigma}{\sigma}.
\end{equation}
Indeed, the integral up to $R$ equals
$[R^\sigma-(R+L)^\sigma+L^\sigma]/\sigma$, and the first two terms tend to zero. Thus
\[
R_a(x)=R_a(a)+\frac1{\Gamma(\sigma)}
\int_{-\infty}^{a}\bigl[(x-t)^{\sigma-1}-(a-t)^{\sigma-1}\bigr]f(t)\dd t
\]
is well-defined, with an absolutely convergent difference integral, and equals the corresponding improper history. For $a\le x<y\le b$ and $M=\|f\|_{\infty;(-\infty,a]}$, the same calculation shows
\[
|R_a(y)-R_a(x)|
\le\frac{M}{\Gamma(1+\sigma)}
\bigl[(y-a)^\sigma-(x-a)^\sigma\bigr].
\]
This proves continuity. Summing over a partition telescopes the right-hand side and yields \eqref{eq:conditionalvariation}.
\end{proof}

\begin{proposition}[Zero-mean periodic Weyl integrals]\label{prop:periodicweyl}
Let $f\in C(\R)$ be $P$-periodic, where $P>0$, and assume
\[
\int_0^P f(t)\dd t=0.
\]
The improper integrals \eqref{eq:weylleft} and \eqref{eq:weylright} exist for every $x$ and $0<\sigma<1$, and are continuous periodic functions. If $f|_{[a,b]}\in C^h([a,b])$ and \eqref{eq:anti} holds on a nondegenerate subinterval of $[a,b]$, then \eqref{eq:extensionbounds} holds for both operators on $[a,b]$.
\end{proposition}
\begin{proof}
For $x\in\R$, the function
\[
A_x(T)=\int_0^T f(x-t)\dd t,\qquad T\ge0,
\]
is $P$-periodic in $T$ and satisfies $|A_x(T)|\le P\|f\|_\infty$. For $U>T>0$, integration by parts yields
\[
\int_T^U t^{\sigma-1}f(x-t)\dd t
=\bigl[t^{\sigma-1}A_x(t)\bigr]_T^U
-(\sigma-1)\int_T^U t^{\sigma-2}A_x(t)\dd t.
\]
Since $\sigma<1$, the term at $U$ tends to zero and the remaining integral converges absolutely. Its tail, including the term at $T$, is bounded by $C_{P,\sigma}\norm f_\infty T^{\sigma-1}$, uniformly in $x$. Near zero, the kernel is integrable and $f$ is bounded. Thus
\[
\Weyl_+^{\sigma}f(x)
=\frac1{\Gamma(\sigma)}\int_0^\infty t^{\sigma-1}f(x-t)\dd t
\]
exists as an improper integral. The truncated integrals are continuous and periodic in $x$, and the uniform tail estimate proves the same properties for their limit.

At $a$, this convergence is precisely the hypothesis of Lemma~\ref{lem:conditionaltail}. Splitting at $a$ yields \eqref{eq:weylsplit} with a continuous history of bounded variation. After affine rescaling, \eqref{eq:extensionbounds} follows from Lemma~\ref{lem:bv} and Theorem~\ref{thm:main}. For the right-sided integral, convergence is proved by replacing $f(x-t)$ by $f(x+t)$, and the dimension bounds are obtained by reflection.
\end{proof}

The zero-mean convention is essential for the ordinary improper periodic integral. A nonzero mean $m$ contributes $m\int_0^R t^{\sigma-1}\dd t=mR^\sigma/\sigma$, while the zero-mean part has a convergent integral by the proof above. Thus the ordinary integral diverges if $m\ne0$. A periodic operator defined by discarding the zero Fourier mode should instead be interpreted as acting on $f-m$; the oscillations and H\"older seminorms are unchanged by this subtraction.

\begin{corollary}[Sharpness in the periodic class]\label{cor:periodicsharp}
For every $0<h\le H<1$, there is a continuous one-periodic zero-mean function $p$, globally $h$-H\"older and uniformly $H$-anti-H\"older on $[0,1]$, for which
\[
\bdim\G(\Weyl_+^{\sigma}p,[0,1])
=\max\{1,2-H-\sigma H/h\},\qquad 0<\sigma<1.
\]
The input graph has dimension $2-H$. The integral has bounded variation on $[0,1]$ when $H+\sigma H/h>1$.
\end{corollary}
\begin{proof}
Periodically extend $f_*$ from Theorem~\ref{thm:sharp}, and let $p$ be that extension minus $m=\int_0^1f_*(t)\dd t$. The extension is continuous because the endpoint values agree. To check its global $h$-H\"older bound, pairs at distance at least one are handled by boundedness. Pairs in one period are handled by the original bound. For pairs across a join and at distance less than one, compare both points to the common endpoint value and use $r_1^h+r_2^h\le2^{1-h}(r_1+r_2)^h$. Subtracting $m$ preserves all oscillations, so the anti-H\"older bound and the input dimension are unchanged.

On $[0,1]$ the local integral is
\[
I_{0+}^{\sigma}p(x)
=I_{0+}^{\sigma}f_*(x)-\frac{m}{\Gamma(1+\sigma)}x^\sigma.
\]
The last term has bounded variation. By Proposition~\ref{prop:periodicweyl} and its proof, $\mathcal W_+^\sigma p$ is obtained by adding a further continuous term of bounded variation. The exact dimension and variation property follow from Lemma~\ref{lem:bv} and Theorem~\ref{thm:sharp}. The construction is independent of $\sigma$.
\end{proof}

\section{Interpretation and scope}\label{sec:scope}

The two exponents have different roles. For $0<r<1$, an upper increment bound $Cr^h$ becomes stronger as $h$ increases, whereas an oscillation lower bound $cr^H$ becomes weaker as $H$ increases. A useful lower bound therefore needs positive H\"older control together with an upper control on the anti-H\"older exponent. A lower bound on the anti-H\"older exponent alone has the wrong direction for this purpose.

The sharp construction displays three regimes. If $H+\sigma H/h<1$, the integral graph has dimension $2-H-\sigma H/h>1$. At equality its graph dimension is one, with the logarithmic covering estimate supplied by \eqref{eq:threecases}. Above the threshold the integral has bounded variation, which is not asserted at the threshold. The Riemann--Liouville, Hadamard, and compactly supported Weyl sharpness examples can all be chosen nonnegative. Thus sign cancellation in those input densities is not necessary for the loss of a nontrivial fractal dimension; the width of the oscillations already supplies the obstruction. The periodic Weyl construction is centered to have zero mean, as required for convergence of its ordinary improper integral.

The pointwise and uniform notions of H\"older regularity should not be conflated. Theorem~\ref{thm:main} assumes the actual two-point estimate \eqref{eq:holder}; a statement about pointwise H\"older exponents alone is not substituted for it. Proposition~\ref{prop:baire} removes uniformity of the constants in a precise local anti-H\"older condition, while retaining a fixed exponent or a uniform upper bound for those exponents. Proposition~\ref{prop:local} explains why a single interval carrying both quantitative assumptions suffices. In each case the quantifiers are part of the hypothesis.

For Hadamard integrals, the compactness and positivity of the observation interval make the logarithmic coordinate change bi-Lipschitz. This proof does not extend automatically to a closed interval with endpoint zero, or to an unbounded interval: convergence and geometry at the omitted endpoint require separate assumptions. For Weyl integrals, a local H\"older hypothesis alone does not control convergence at infinity. Conditions \eqref{eq:tailintroleft} and \eqref{eq:tailintroright} are explicit sufficient conditions, while the periodic result uses cancellation instead. They should not be conflated. The bounds apply to the graph restricted to the stated compact observation interval.

Finally, these results concern mesh-based graph dimensions. Their proofs rely on sums of oscillations and on invariance under addition of continuous bounded-variation functions. They do not yield the analogous lower bound for Hausdorff dimension. The exact extremal formula obtained here is
\[
\inf_f\lbdim\G(I_{0+}^{\sigma}f,[0,1])
=\max\{1,2-H-\sigma H/h\},
\]
where the infimum is taken over functions satisfying the hypotheses of Theorem~\ref{thm:main} for fixed $h,H$ and arbitrary positive admissible constants. It is attained by the explicit family in Section~\ref{sec:sharp}. In this sense the unequal-exponent lower bound is optimal. The same extremal value is attained in the Hadamard, compactly supported Weyl, and zero-mean periodic Weyl classes by Propositions~\ref{prop:hadsharp} and \ref{prop:weylsharp} and Corollary~\ref{cor:periodicsharp}, respectively. These extensions follow from the Riemann--Liouville theorem. A separate bibliographic comparison would be needed to establish publication priority.

\pdfbookmark[2]{Comparison with other graph-dimension results}{comparison}
\subsection*{Comparison with other graph-dimension results}
For a continuous function whose graph has upper box dimension $s>1$, Wu, Ji, and Zhang \cite{WuJiZhang} obtained the upper bound $s-(s-1)\sigma$ for the Riemann--Liouville integral of order $0<\sigma<1$. In Theorem~\ref{thm:main} the upper bound is determined by the increment exponent $h$. The lower bound also uses the oscillation exponent $H$ on a single subinterval. When $h=H$, Corollary~\ref{cor:equal} recovers the linear law addressed in \cite[Remark~3.13]{LiangSuGeometric}. For unequal exponents the sharp lower bound contains the factor $H/h$.

Exact box dimensions have also been studied for fractional integrals of Weierstrass functions by Zhou, Yao, and Su \cite{ZhouYaoSu}, and for fractional integrals and derivatives of Besicovitch functions by He and Zhou \cite{HeZhou}. For linear fractal interpolation functions, Ruan, Su, and Yao \cite{RuanSuYao} obtained a dimension law in cases where the integral retains the interpolation structure. The construction in Theorem~\ref{thm:sharp} serves a different purpose: for each pair $(h,H)$ it attains the lower bound for every integration order using the same input. This permits optimality to be tested over the entire class in Theorem~\ref{thm:main}.

The dimension being estimated is important in interpreting the Hadamard extension. Yao and coauthors \cite{YaoHadamard} studied Hausdorff dimension for the Hadamard integral of a particular fractal function. Theorem~\ref{thm:extensions} concerns lower and upper box dimensions, both of which are preserved by the logarithmic change of coordinates on the compact observation interval.

\clearpage
\pdfbookmark[1]{References}{references}

\end{document}